\documentstyle[12pt]{amsart}

\makeatletter
\def\verbatim{\interlinepenalty\@M \@verbatim
  \leftskip\@totalleftmargin\advance\leftskip2pc
  \frenchspacing\@vobeyspaces \@xverbatim}
\makeatother
\newtheorem{theorem}{Theorem}[section]

\theoremstyle{definition}

\numberwithin{equation}{section}

\title[Isometric shifts]{Isometric shifts and metric spaces}

\author{Jes\'us Araujo}
\address{Departamento de Matem\'aticas,
Estad\'{\i}stica y Computaci\'on\\ Universidad de Cantabria\\
Facultad de Ciencias\\ Avda.
de los Castros, s. n.\\ E-39071 Santander, Spain}

\email{araujoj@@unican.es}
\thanks{2000 {\em Mathematics Subject Classification}.
Primary 47B38; Secondary 54D65, 46J10.}
\thanks{Research of the first author was partially supported by
the Spanish Direcci\'on General de Investigaci\'on Cient\'{\i}fica
y T\'ecnica (DGICYT, PB98-1102).}

\author{Juan J. Font}
\address{Departamento de Matem\'aticas\\
Universitat Jaume I\\ Campus Penyeta Roja\\ E-12071 Castell\'on, Spain}
\email{font@@mat.uji.es}
\thanks{Research of the second author was partially supported by
Fundaci\'o Caixa Castell\'o.}

\date{}                                          
\begin{document}

\maketitle

\maketitle







\begin{abstract}
Let $M$ be a complete metric space. If $C^*(M)$ admits an isometric shift,
then $M$ is separable.
\end{abstract}

\section{Introduction}

Shift operators play an important role in many disciplines 
such as Perturbation Theory, Engineering Mathematics, Scattering Theory,
Stochastic Processes, \ldots (see \cite{N}).
R.M. Crownover (\cite{Cr}) was the first to extend the definition
of shift operator from separable Hilbert spaces
to arbitrary Banach spaces without using
basis. Namely, if ${\cal K}$ is a Banach space, then 
$T:{\cal K} \longrightarrow {\cal K}$ is said to be an {\em (isometric) 
shift operator} if 
\begin{enumerate}
\item $T$ is a linear isometry,
\item The codimension of $T({\cal K})$ in ${\cal K}$ is $1$,
\item $\bigcap ^{\infty}_{n=1} T^{n}({\cal K})=\{0\}.$
\end{enumerate}
If Condition $3$ is removed, then we have a {\em codimension 1
linear isometry}.

In \cite{GHJR}, Gutek, Hart, Jamison and Rajagopalan
extended many of the results obtained by J.R. Holub in \cite{Hl}
concerning isometric shift operators on the Banach
space $C(X)$ 
($X$ compact Hausdorff). First, they classified codimension 1 linear
isometries on $C(X)$ using the following result: 
let $T:C(X) \longrightarrow C(X)$ be a codimension 1 linear isometry. 
Then there exists a closed subset $X_0$ of $X$ such that either 

(i) $X_0=X\setminus \{p\}$ 

\noindent
where $p$ is an isolated point of $X$, or 
\indent

(ii) $X_0=X$

\noindent
and such that there exists a continuous map $h$ of $X_0$ onto $X$ and a 
function $a\in C(X_0)$, $\left|a \right| \equiv 1$, such that 
\[(Tf)(x)=a(x)\cdot f(h(x))\]
for all $x\in X_0$.
\indent

The proof of this result is based on a well known theorem of Holszty\'{n}ski 
(\cite{Hz}). Those isometries that satisfy Condition (i) are said to be of 
type I. Those satisfying Condition (ii) are said to be of type II. 
These two classes are not disjoint. 
Farid and Varadarajan (\cite{FV}) devoted part of their paper to clarify
the above classification. Finally, in \cite{F} the author
proposes an alternative (disjoint) classification based on the
separation properties of the range of $T$. Thus, 
$T$ is of type II if and only if $T(C(X))$ 
separates  all the points of $X$ except two and is of type I
which is not of type II if and only if $T(C(X))$
separates  all the points of $X$. 

Codimension 1 linear isometries on arbitrary function algebras
have also been studied and classified in \cite{AF2} by using the
results in \cite{AF1}.
Recently, K. Izuchi \cite{Iz} has characterized Douglas
algebras which admit codimension 1 linear isometries, thus
solving the conjecture settled in \cite{AF2}.

Another question which has also been addressed in the context
of isometric shifts is the
characterization of those compact Hausdorff spaces $X$ which admit
such operators, that is, the existence of
isometric shifts on $C(X)$.
In \cite{GHJR}, the authors
proved that nonseparable spaces without isolated points do not
admit isometric shifts and even that there is no nonseparable space
which admits isometric shifts of type II. They left open the question
of the existence of spaces without isolated points which admit
isometric shifts. This question was answered in the
positive by R. Haydon (\cite{Ha}). He proved the existence 
of isometric shifts of
type II when $X$ is either connected or the Cantor set.
However it is still an open question whether there exists a nonseparable
compact space $X$ which admits an isometric shift.
In this paper we show that no nonseparable metric
(noncompact) space admits isometric shifts. We also provide an example
of an isometric shift (of type I) with several interesting features.

\section{Preliminaries}

Let $\Bbb{K}$ denote the field of
real or complex numbers. 
If $X$ is a compact (respectively locally compact) Hausdorff space, then
$C(X)$ (respectively $C_0(X)$)
stands for the Banach space of all 
$\Bbb{K}$-valued continuous functions defined on $X$ (respectively which vanish at infinity), 
equipped with its usual supremum norm.
If $M$ is a metric space, then we shall write $C^*(M)$ to denote the normed
space of all {\em bounded} $\Bbb{K}$-valued continuous functions defined on $M$.
As usual, $\beta M$ stands for the Stone-\v{C}ech compactification of $M$.
Given $f\in C(X)$, we shall consider that $c(f)$ is 
its cozero set.

If $U$ is a subset of $X$, then ${\rm cl}_X(U)$ and ${\rm int}_X(U)$ denote
its closure and its interior in $X$, respectively.

\section{Isometric shifts on $C^*(M)$ and separability}

\bigskip
Let $M$ be a complete metric space and let $T:C^*(M)\longrightarrow C^*(M)$
be an isometric shift. Then $T$ induces an isometric shift
(which we continue to denote by $T$) on $C(\beta M)$.

\begin{theorem}  \label{main}
Let $M$ be a complete metric space. If $C^*(M)$ admits an isometric
shift $T$, then $M$ is separable.
\end{theorem}

{\em Proof.}
Let us first assume $T$ to be of type II which is not of type I.
According to \cite[Lemma 2.2]{GHJR}, the map $h: \beta M \longrightarrow
\beta M$ is a surjective continuous map such that there exists $x_0 \in
\beta M$ in such a way that $h^{-1} (\{x\})$ consists of just one point
for every $x \in \beta M \setminus  \{x_0 \}$ and $h^{-1} (\{x_0 \})$
consists of two points of $ \beta M$, say $x_1 ,x_2$.
Also, since
$T$ is not of type I, then the points $x_1$ and $x_2$ are not isolated.
Furthermore, in \cite[Theorem 2.5]{GHJR}, it is proven that the set
\[{\cal D} := \bigcup_{k=-\infty}^{\infty} h^k (\{x_0\})\] is a countable
dense subset in $\beta M$. We are going to see that this set is contained in
$M$ and in this way we give an explicit countable dense subset in $M$.

First we have, by \cite[Theorem 2.6]{GHJR}, that if $\beta M/ R$ is
the quotient space for the equivalence relation defined as $x R y$
whenever $h(x) =h(y)$, then the map $h^R : \beta M/ R \longrightarrow
\beta M$ sending each class $x^R$ into the image $h(x)$ of any $x \in x^R$
is a surjective homeomorphism. This implies in particular that the image
of a $G_{\delta}$-point in $\beta M /R$ is a $G_{\delta}$-point in $\beta M$
and vice versa. Let us recall that the only points
in $\beta M$ which are $G_{\delta}$ are those in $M$.

Let us check which the $G_{\delta}$-points in $\beta M /R$ are. Suppose that
$x^R \in \beta M / R$ satisfies that there exists $x \in M$ with
$x  \in x^R$. Clearly, if $x^R$ is the singleton $\{x\}$,
then $x^R$ is $G_{\delta}$. Otherwise, as we remark above,
$x^R$ consists of two points, $x_1, x_2$, and is the only point in
$\beta M /R$ which is not a singleton.
Then it is apparent that $x^R$ is $G_{\delta}$
if and only if both $x_1,x_2 \in M$.

Suppose next that
$x^R = \{x_1, x_2\}$, and that $x_1 \in M$.
Since $T$ is not of type I, then the points $x_1$ and $x_2$ are clearly not
isolated. Thus, there exists a sequence $(y_n)$
in $M \setminus \{x_1, x_2\}$ converging to $x_1$. Also each $y_n$ is
a $G_{\delta}$-point, and consequently so is $h^R (y_n)=h(y_n)$,
that is, the sequence $(h(y_n))$ is contained in $M$, and converges
to $x_0$. But this implies in particular that $x_0 \in M$ (\cite[Theorem 8.3.2]{Wil}).
Conversely, if we assume that $x_0\in M$, then $x_0$ is a
$G_{\delta}$-point of $\beta M$ and, consequently,
so is $h^{-1} (\{x_0\})=\{x_1, x_2\}$. This implies, as stated above,
that both $x_1$ and $x_2$ belong to $M$.
Summing up, we have proven that $x_0 \in M$ if and only if
$x_1 \in M$ or $x_2 \in M$, and that this fact yields
$x_1 , x_2 \in M$.

Let us now assume that $x_0 \notin M$, which is to say that
$x_1 , x_2 \notin M$. Then it is easy to check
that the restriction of the map $h$ to $M$, $h:M \longrightarrow M$,
is bijective and continuous, and its inverse $h^{-1}:M \longrightarrow M$
is also continuous. Consequently, the map $T : C^* (M) \longrightarrow
C^* (M)$ sending each $f$ into $a \cdot f \circ h$, $\left| a \right|
\equiv 1$, is clearly a surjective linear isometry, that is, it is
{\em not} a codimension $1$ isometry, against our hypothesis. We
deduce that $x_0$ must belong to $M$, and consequently $x_1 , x_2$
belong to $M$. Hence, $h^{-1} (\{x_0 \}) \subset M$.

A similar reasoning leads to the fact that $h^k (\{x_0\}) \subset M$
for every integer $k$. That is, ${\cal D} \subset M$, as was
to be proved.

\bigskip
Let us now assume that $T:C(\beta M) \longrightarrow C(\beta M)$
is of type I. Thus, there exist an isolated
point $p\in \beta M$ and a homeomorphism (\cite[Lemma 2.2]{GHJR}) $h$ of $\beta M\setminus \{p\}$ onto
$\beta M$ and a
function $a\in C(\beta M\setminus \{p\})$,
$\left|a \right| \equiv 1$, such that
\[(Tf)(x)=a(x)\cdot f(h(x)) \]
for all $x\in \beta M\setminus \{p\}$.
Consider the set $A=\{p,h^{-1}(p),h^{-2}(p),...\}$.
Then $Y:=\beta M \setminus {\rm cl}_{\beta M}(A)$  is a locally compact space
and $h:Y\longrightarrow Y$ is a surjective homeomorphism.
Hence we have a surjective isometry
$S:C_0(Y)\longrightarrow C_0(Y)$ defined to be 
\[(Sf)(x)=\hat{a}(x)\cdot f(h(x)),\]
where $\hat{a}$ is the restriction to $Y$ of $a$. 

For any $f\in C_0(Y)$, we can define a function
$\hat{f}\in C(\beta M)$ such that
$\hat{f}=f$ on $Y$ and $0$ on $\beta M\setminus {Y}$.
As a consequence, a linear continuous functional $\mu$ (indeed a regular complex measure)
can be defined on $C_0(Y)$ to be $\mu(f):=(T\hat{f})(p)$.

\bigskip
{\bf Claim 1.} {\em Assume that there is $f\in C_0(Y)$
such that $\mu(f)=0$ and $(\mu \circ S^{-n})(f)=0$ for all $n\in \Bbb{N}$.
Then $f\equiv 0$.}

\bigskip
Let us suppose, contrary to what we claim, that there
is $f\in C_0(Y)$, $f\neq 0$, 
such that $\mu(f)=0$ and $(\mu \circ S^{-n})(f)=0$ for all $n\in \Bbb{N}$.
Let us check that $\hat{f}\in R(T^n)$ for all $n\in \Bbb{N}$.

Since $S:C_0(Y)\longrightarrow C_0(Y)$ is a surjective isometry,
there is $g\in C_0(Y)$ such that $S(g)=f$.
If $x\in Y$, then
\[(T\hat{g})(x)=a(x)\cdot \hat{g}(h(x))=
\hat{a}(x)\cdot g(h(x))=(Sg)(x)=f(x)=\hat{f}(x).\]
That is, $T\hat{g}=\hat{f}$ on $Y$.

On the other hand,
$(T\hat{g})(p):=\mu(g)=\mu(S^{-1}f)=(\mu\circ S^{-1})(f)$.
By assumption, $(\mu\circ S^{-1})(f)=0$. Hence, $(T\hat{g})(p)=0=\hat{f}(p)$.

Next, from the representation of the isometric shift $T$, we know that
$(T\hat{g})(h^{-n}(p))=a(h^{-n}(p))\cdot \hat{g}(h^{-n+1}(p))$, but
$h^{-n+1}(p)\in \beta M\setminus Y$, which is to say that
$\hat{g}(h^{-n+1}(p))=0$.

Finally, it is apparent, from the above two paragraphs and from density,
that $T\hat{g}\equiv 0$ on $\beta M\setminus Y$. Hence, gathering the
information above, we infer that $T\hat{g}=\hat{f}$, i.e.,
$\hat{f}\in R(T)$.

Let us next check that $\hat{f}\in R(T^2)$. To see this, it suffices
to prove that $\hat{g}\in R(T)$. Since $S$ is surjective, 
there is $g_1\in C_0(Y)$ such that $S(g_1)=g$. Furthermore
$(T\hat{g_1})(p):=\mu(g_1)=\mu(S^{-2}f)=(\mu\circ S^{-2})(f)=0=\hat{g}(p)$.
Hence, as above, we deduce that $T\hat{g_1}=\hat{g}$.
In like manner, we can obtain $g_2, g_3, \ldots, g_n, \ldots$ to show that
$\hat{f}\in R(T^n)$ for all $n\in \Bbb{N}$.  This fact contradicts the
definition of isometric shift and the proof of Claim 1 is complete.

\bigskip
It is well-known that every regular complex measure $\theta$ can be written
as $\theta=(\theta_1-\theta_2)+i(\theta_3-\theta_4)$, where
$\theta_i$, $i=1,2,3,4$,
are regular positive measures.
Hence each of the regular complex measures
$\mu, \mu \circ S^{-1}, \mu \circ S^{-2}, \ldots , \mu \circ S^{-n},
\ldots$ can be divided into four regular positive measures.
As a consequence, we get a new sequence of regular positive measures,
which we shall denote by $\{\mu_n\}_{n\in \Bbb{N}}$. With no loss
of generality, we can assume that all these measures are normalized.

Since the space of regular measures on a locally compact space is a Banach space,
we can define
a regular positive measure as follows:
\[\eta:=\sum_{n=1}^{\infty}\frac{\mu_n}{2^n}.\]

\bigskip
{\bf Claim 2.} {\em For every nonempty open subset $U$ of $Y$, $\eta(U)>0$.}

\bigskip
Let us suppose that there exists a nonempty open subset $U$ of $Y$
such that $\eta(U)=0$. Hence we can find $f\in C_0(Y)$, $f\neq 0$,
such that $c(f)\subset U$. Consequently,
\[\mu_n(f)=\int_{Y}fd\mu_n=0\]
for all $n\in \Bbb{N}$. Finally, Claim 1 yields $f\equiv 0$, a contradiction.

\bigskip
Let us now define a (open) subset $N:=M\setminus
{\rm cl}_{\beta M}(\{p,h^{-1}(p),h^{-2}(p),...\})$ of $M$.
Next we consider the family, say 
${\cal F}_1$, of all subsets $B$ of $N$
which satisfy the following property:
if $x,y\in B$, then $d(x,y)\ge 1$ or $d(x,y)=0$, where $d$ denotes
the metric in $N$ induced from $M$.
Let us choose a chain $(A_{\alpha})_{\alpha}$ of elements of
${\cal F}_1$  ordered by inclusion. Since
\[\bigcup_{\alpha} A_{\alpha}\in {\cal F}_1,\]
Zorn's lemma yields a maximal element, say $M_1$.

\bigskip
{\bf Claim 3.} {\em $M_1$ is a countable set.}

\bigskip
Assume the contrary. Then there exists an uncountable family
$\triangle$ of indexes such that
\[M_1=\{x_{\alpha}:\alpha\in \triangle\}.\]
Since $N$ is an open subset of $M$, there is, for each
$\alpha\in \triangle$, a constant $M_{\alpha}>0$ such that the open
ball $B(x_{\alpha},M_{\alpha})\subset N$.

Next, for each $\alpha\in \triangle$, take
\[m_{\alpha}:=\inf\left\{\frac{1}{3},M_{\alpha}\right\}\]
and consider the open ball $B(x_{\alpha},m_{\alpha})$.
It is clear, from the definition of ${\cal F}_1$, that if
$\alpha,\beta\in \triangle$, $\alpha\neq \beta$, then
\[B(x_{\alpha},m_{\alpha})\cap B(x_{\beta},m_{\beta})=\emptyset.\]
Now, for every $\alpha\in \triangle$, we can define the set
\[V_{\alpha}:={\rm int}_{\beta M}({\rm cl}_{\beta M}(B(x_{\alpha},m_{\alpha})).\]
It is apparent that $V_{\alpha}\cap M=B(x_{\alpha},m_{\alpha})$
for each $\alpha\in \triangle$
and that $V_{\alpha}\cap V_{\beta}=\emptyset$ if
$\alpha\neq \beta$. Furthermore,
each $V_{\alpha}$ is contained in $Y$ since $Y$ is open.


Summarizing, we have found an uncountable pairwise disjoint family
of open subsets $\{V_{\alpha}:\alpha\in \triangle\}$ in $Y$.

We know, by Claim 2, that $\eta(V_{\alpha})>0$ for all
$\alpha\in \triangle$. Hence, there is $n_0\in \Bbb{N}$ such that
the set
\[\gamma:=\left\{\alpha\in \triangle: \eta(V_{\alpha})>\frac{1}{n_0}\right\}\]
is not countable since neither is $\triangle$. 
Let us choose a countable subset
$\{\alpha_1,\alpha_2,...,\alpha_n,...\}$ of indexes in $\gamma$.
Then,
\[\eta(\bigcup_{n=1}^{\infty}V_{\alpha_n})=
\sum_{n=1}^{\infty}\eta(V_{\alpha_n})=+\infty.\]
This contradiction completes the proof of Claim 3.

\bigskip
As in the paragraph before Claim 3, we can define, for every
$n\in \Bbb{N}$, the family ${\cal F}_n$ of all subsets $B$ of $N$
which satisfy the following property:
if $x,y\in B$, then $d(x,y)\ge 1/n$ or $d(x,y)=0$.
In like manner, we obtain, for every $n\in \Bbb{N}$, a maximal element
$M_n$ of ${\cal F}_n$ which turns out to be countable.

Let us now see that the countable set
\[{\cal D}:=\bigcup_{n=1}^{\infty}M_n\]
is dense in $N$. To this end, choose $x\in N\setminus {\cal D}$ and
$\epsilon>0$.
Then there exists
$m_0\in \Bbb{N}$ such that $\frac{1}{m_0}<\epsilon$. Since
$x\notin {\cal D}$, then $x\notin M_{m_0}$. This facts implies
the existence of $y\in M_{m_0}$ such that $d(x,y)< 1/{m_0}$. That is,
there is an element $y$ of ${\cal D}$ in the open ball
$B(x,\epsilon)$ and the density of ${\cal D}$ in $N$ follows.

Finally, it is clear that the countable set
\[{\cal D}\cup \{p,h^{-1}(p),h^{-2}(p),...\}\]
is dense in $M$ and we are done.


\section{Example}

\bigskip
In \cite{Ha}, Haydon showed a method to provide isometric shifts of type II.
However, it is remarkable the scarcity of examples of isometric shifts of
type I.
In this final section we provide an example of an isometric shift
of type I, which is not of type II, with several additional features.
Indeed, in \cite{GHJR}, the authors raised
the question whether, for an isometric shift of type I,
the set $D:=\{p,h^{-1}(p),h^{-2}(p),...\}$ was always
dense in $X$. The question
was answered in the negative by Farid and Varadarajan (\cite{FV})
by providing an example of an isometric shift of type I such that
$X\setminus {\rm cl}_{X}(D)$ was a (finite) nonempty subset.
Our example shows somehow that $D$ can be far from being dense in $X$
in the sense that $X\setminus {\rm cl}_{X}(D)$ is uncountable. Our $X$ also has,
contrary to what Holub conjectured in \cite{Hl}, an infinite
connected component (see also \cite[Corollary 2.1]{GHJR}).

\bigskip
{\bf Example.} Let $\partial D$ denote the unit circle in ${\Bbb C}$, and
let
\[X = \partial D \cup \left\{ \frac{1}{n}: n \in {\Bbb N}, n  \ge 2\right\}
\cup \left\{0\right\}.\]
It is clear that $X$ is a compact metric space. Let us show that $X$ admits
an isometric shift of type I by constructing it explicitely.

Let $T: C(X) \rightarrow C(X)$ be the following operator.
Take any $f \in C(X)$ and define, for each $e^{i \theta} \in \partial D$,
\[(Tf) (e^{i \theta}):=f (e^{i (\theta + \sqrt{2})}).\]

It is clear that, given any $e^{i \theta} \in \partial D$,
the sequences $(e^{i (\theta + 2n \sqrt{2})})$  and
$(e^{i (\theta + (2n -1) \sqrt{2})})$ are dense in
$\partial D$. Then we take in $\partial D$ the point $1= e^{i0}$.

Clearly the evaluation map $\delta_1$ is continuous in $C(X)$ and
its norm is equal to $1$. So, for $f \in C(X)$, we define
\[(Tf) (1/2) := - (\delta_1 + \delta_{e^{-i \sqrt{2}}}) (f)/2
=- f(1)/2 -f(e^{-i \sqrt{2}})/2.\]
Next, for $n \ge 3$, we
define
\[(Tf) (1/n) := - f(1/n-1),\]
and
\[(Tf)(0):= - f(0).\]

It is clear that $Tf \in C(X)$ and that $T$ is an isometry.
In fact $T$ is a codimension 1 linear isometry of type
$I$ (being $p= 1/2$), which is not of type II since the range of $T$
separates all the points of $X$ (see \cite{F}).
Let us see that it is also a shift operator.

Suppose that $g \in C(X)$ satisfies
$g \in \bigcap_{n=1}^{\infty} R(T^n)$. We have to prove that $g=0$.
First we have that $(g(1/n))$ must be a convergent sequence, and
it converges to the value $g(0)$. Also, we have that
$g(1/2) = - (T^{-1} g) (1)/2 - (T^{-1} g) (e^{- i \sqrt{2}})/2$,
by construction. In the same way $g (1/3) = - (T^{-1} g) (1/2) =
(T^{-2} g) (1) /2 + (T^{2} g) (e^{-i \sqrt{2}}) /2=
(T^{-1} g) (e^{-i \sqrt{2}})/2 + (T^{-1} g) (e^{-i 2 \sqrt{2}})/2$
and, in general, for $n \ge 2$, $n \in {\Bbb N}$, \[g(1/n) =
(-1)^{n+1}( (T^{-1} g) (e^{- (n-1) i \sqrt{2}})/2 +
(T^{-1} g) (e^{- (n-2) i \sqrt{2}})/2 ).\]

In particular, we have that the sequence 
\[\left((-1)^{n+1} \frac{(T^{-1} g) (e^{- (n-1) i \sqrt{2}}) +
(T^{-1} g) (e^{- (n-2) i \sqrt{2}})}{2}\right)\] must converge to $g(0)$,
because $g$ is continuous.

On the other hand, by the density of points of the form
$e^{i 2n \sqrt{2}}$, $n \in {\Bbb N}$, we have that given any point
$z_0 \in \partial D$, there exists a sequence $(n_k)$ of even numbers
such that $(e^{- i  n_k \sqrt{2}})$ converges to $z_0$ as $k$ tends
to infinity. Also
$(e^{- i  (n_k -1) \sqrt{2}})$ converges to $z_0 e^{-i \sqrt{2}}$.
Since $T^{-1} g$ is continuous, this implies that $((T^{-1} g)
(e^{- i  n_k \sqrt{2}}))$  goes to $(T^{-1} g) (z_0)$, and that
$((T^{-1} g) (e^{- i  (n_k -1) \sqrt{2}}))$  goes to $(T^{-1} g)
(z_0 e^{-i \sqrt{2}})$. We deduce that
\[\left(\frac{ (T^{-1} g)
(e^{- i  (n_k -1) \sqrt{2}}) + (T^{-1} g) (e^{- i  n_k  \sqrt{2}})}{2}\right)\]
converges to
\[\frac{(T^{-1} g) (z_0 ) + (T^{-1} g) (z_0 e^{-i \sqrt{2}})}{2}.\] 
On the other hand, we know that the above sequence converges to $g(0)$.
But a similar approach can be taken for a sequence of odd natural numbers
$(m_k)$ instead of $(n_k)$. In this case we will obtain that
\[\left(\frac{ (T^{-1} g) (e^{- i  (m_k -1) \sqrt{2}}) + (T^{-1} g)
(e^{- i  m_k  \sqrt{2}})}{2}\right)\] converges to
\[\frac{(T^{-1} g) (z_0 ) + (T^{-1} g) (z_0 e^{-i \sqrt{2}})}{2},\]
and on the other hand, it must converge to $-g(0)$. As a consequence,
we deduce that $g(0) =-g(0) =0$, and that, for every $z_0 \in \partial D$,
\[(T^{-1} g) (z_0) + (T^{-1} g) (z_0 e^{-i \sqrt{2}}) =0 .\] In particular,
this implies that for every $z_0 \in \partial D$, $(T^{-1} g)
(z_0 e^{-i \sqrt{2}}) = (T^{-1} g) (z_0 e^{i \sqrt{2}})$. Consequently,
the sequence \[\left((T^{-1} g) (e^{i 2n \sqrt{2}})\right)\] is constant. By the
density of points $e^{i 2n \sqrt{2}}$, $n \in {\Bbb N}$, we conclude
that $T^{-1} g$ is constant on $\partial D$. In particular, this implies
that the sequence $(\left|g(1/n)\right|)$ is constant.
Since it converges to $\left|g(0)\right|=0$,
we conclude that $g(1/n)=0$ for every $n \in {\Bbb N}$. As a consequence
it is easy to see that $T^{-1} g \equiv 0$ on $\partial D$. But this
clearly implies that $g=0$, as we wanted to prove.

\end{document}